\documentclass[12pt]{amsart}
\usepackage{times,amsmath,amsfonts,amssymb,amscd,a4wide}
 
\usepackage{verbatim}
\usepackage{hyperref}
\usepackage[all,cmtip]{xy}
 
\newtheorem{theor}{Theorem}

\newtheorem{theo}{Theorem}

\newcommand{\Sph}{\mathbb{S}}

\newcommand{\Z}{\mathbb{Z}}
\newcommand{\R}{\mathbb{R}}
\newcommand{\C}{\mathbb{C}}
\newcommand{\F}{\mathbb{F}}
\newcommand{\NN}{\mathbb{N}}
\newcommand{\GL}{\mathbf G\mathbf L}

\newcommand{\Aut}{\operatorname{Aut}}

\newcommand{\Sym}{\operatorname{Sym}}
\newcommand{\Tr}{\operatorname{Tr}}

\newcommand{\lb}{\left<}
\newcommand{\rb}{\right>}

\renewcommand{\Re}{\mathrm{Re}}
\renewcommand{\Im}{\mathrm{Im}}

\newcommand{\eg}{{\it e.g. }}
\newcommand{\ie}{{\it i.e. }}

\newtheorem{theorem}{Theorem}[section]
\newtheorem{proposition}[theorem]{Proposition}
\newtheorem{remark}[theorem]{Remark}
\newtheorem{corollary}[theorem]{Corollary}

\author{Renaud Coulangeon}

\address{Institut de Math\'ematiques de Bordeaux,
Universit\'e Bordeaux I, 351, cours de la Lib\'eration \\
33405 Talence, France}

\title{Spherical designs and zeta functions of lattices.}
\begin{document}

\begin{abstract}
We set up a connection between the theory of spherical designs and the question of minima of Epstein's zeta function. More precisely, we prove that a Euclidean lattice, all layers of which hold a $4$-design, achieves a local minimum of the Epstein's zeta function, at least at any real $s>\frac{n}{2}$. We deduce from this a new proof of Sarnak and Str\"ombergsson's theorem asserting that the root lattices $\mathbb{D}_4$ and $\mathbb{E}_8$, as well as the Leech lattice $\Lambda_{24}$, achieve a strict local minimum of the Epstein's zeta function at any $s>0$. Furthermore, our criterion enables us to extend their theorem to all the so-called \textit{extremal modular lattices} (up to certain restrictions) using a theorem of Bachoc and Venkov, and to other classical families of lattices (\eg the Barnes-Wall lattices).
\end{abstract}
 
\maketitle
\section*{Introduction.} 
A spherical design is a finite set of points on a sphere which is
\textit{well-distributed}, in the sense that it allows numerical
integration of functions on the sphere up to a certain accurracy,
the so-called \textit{strength} of the design. More precisely, $X
\subset \Sph^{n-1}$ is a $t$-design if for all homogeneous polynomial
of degree $\leq t$
\begin{equation}\int_{\Sph^{n-1}}f(x)dx=\frac{1}{|X|}\sum _{x \in X}f(x)
\end{equation}
where $\Sph^{n-1}$ stands for the unit sphere in $\R^n$ endowed with
its canonical measure $dx$, normalized so as $\int_{\Sph^{n-1}}dx=1$.

One classical way to build such designs is to consider the set of
vectors of a given length in a Euclidean lattice $L$ in $\R^n$,
rescaled so as to lie in $\Sph^{n-1}$. It has long been observed that there is a link
between the ability of getting designs of high strength in that way
and classical properties of the underlying lattice e.g. density,
symmetries, theta series. In this connection one can quote Venkov's
remarkable theorem asserting that if the set of \textit{minimal
vectors} (non zero vectors of minimal length) of a lattice is a
$4$-design, then the lattice is extreme in Vorono\"i's sense \ie it
achieves a local maximum of the packing density function \cite{V}.

There are many examples of lattices for which not only the set of minimal vectors holds a design, but all sets of vectors of any given length actually do. This happens for instance with the so-called \textit{extremal
modular lattices}, as shown by Bachoc and Venkov using theta series with spherical coefficients (see\cite{BV} and Section \ref{ex} below). Another instance of this phenomenon is when the automorphism group of the lattice is "big enough" to satisfy Goethals and Seidel's theorem (see \cite[Th\'eor\`eme 6.1.]{GS}), or Section \ref{ex} below). In all these cases, since not only minimal vectors, but all layers are involved, one would expect some more consequences on the associated packing than local optimality for the density. One aim of the present paper is to provide an interpretation of this phenomenon in terms of Epstein's $\zeta$-function.
 
The Epstein zeta function of a lattice $L$ is defined, for $s\in\C$ with $\Re(s)
> \frac{n}{2}$, as
\begin{equation}
\zeta(L,s):= \sum_{x \in L - \{0\}} ||x||^{-2s}
\end{equation} and admits a meromorphic
continuation to the complex plane with a simple pole at
$s=\frac{n}{2}$. The question as to which $L$, for fixed $s >0$ ($s
\neq \frac{n}{2}$) minimizes $\zeta(L,s)$, has a long history,
starting with Sobolev's work on numerical integration (\cite{So}), and a series
of subsequent papers by Delone and Ryshkov among others. Of course,
this question makes sense only if one restricts to lattices of fixed
covolume, say $1$, since
\begin{equation}
\forall \lambda > 0 \ \ \zeta(\lambda L,s)= \lambda^{-2s}\zeta(L,s).
\end{equation}
From now on, we denote $\mathcal L_n ^{\circ}$ the set of lattices of
determinant (covolume) $1$. Another, undoubtedly  more important, reason to investigate this question, is its connection with Riemannian geometry : if $X$ is a compact Riemannian manifold, its height $h(X)$ is defined as $\zeta_{X}'(0)$, where $\zeta_X(s)$ is the so-called \textit{zeta regularisation} of the determinant of the Laplacian, \ie $\zeta_{X}(s)=\sum_{\lambda_j \neq 0}\lambda_j ^{-s}$, where the $\lambda_j$ are the eigenvalues of the Laplacian on $X$. When $X$ is a flat torus $\R^n/L$, with $L$ a full rank lattice in $\R^n$, then $\zeta_{X}(s)$ is the same as $\zeta(L^{*},s)$, up to a constant (see \cite{Ch} and \cite{Sa}). In this context, a natural question is to find lattices achieving a minimum of this height function restricted to flat tori (the existence of such a minimum is shown in \cite{Ch}), which amounts to minimize  $\zeta'(L^{*},0)$ for $L \in \mathcal L_n ^{\circ}$.

We say that a lattice $L_0 \in \mathcal
L ^{\circ}$ is $\zeta$-\textit{extreme} at $s \in \R$ if it achieves a strict local minimum of the function $L \mapsto \zeta(L,s)$, $L
\in \mathcal L_n ^{\circ}$. Delone and Ryshkov's obtained a
characterization of lattices in $\mathcal L_n ^{\circ}$ that are $\zeta$-\textit{extreme} at $s$ for any large
enough $s$ (\cite[Theorem 4]{DR}). One of the conditions is that all layers of the lattice
hold a $2$-design (see Section \ref{basics}).
 
One would naturally ask for an explicit value $s_0$ such that $L_0$
is $\zeta$-\textit{extreme} at $s$ for any $s \geq
s_0$, but unfortunately, it is not possible to derive such an $s_0$
from Delone and Ryshkov's theorem, nor from its proof . However, using more
sophisticated tools, Sarnak and Str\"ombergsson proved in \cite{Sa-Sa}
the following theorem
\begin{theo}[Sarnak \& Str\"ombergsson (2006)] The $\mathbb{D}_4$ lattice (rescaled so as to have determinant $1$), the $\mathbb{E}_8$-lattice and the Leech Lattice $\Lambda_{24}$
are $\zeta$-\textit{extreme} at $s$ for any $s > 0$, and the associated  tori achieve a strict local minimum of the height function on the set of flat tori of covolume $1$ and dimension $4$, $8$ and $24$ respectively.
 
\end{theo}
Their proof relies essentially on a certain property of the
automorphism group of those lattices which is shown to imply the
desired property, at least for $s>\frac{n}{2}$ (the proof for $s$ in
the "critical strip" $0<s<\frac{n}{2}$ is more involved and requires
some extra arguments). Inspired by Delone and Ryshkov's theorem, one may ask
for an explanation of this result in terms of spherical designs.
This is precisely the aim of our main theorem, which we now state
\begin{theor}\label{mt}
Let $L \in \mathcal L_n ^{\circ}$ be such that all its \textrm{layers}
hold a $4$-design. Then $L$ is $\zeta$-extreme at $s$ for
any $s >\frac{n}{2}$, and the torus associated with its dual $L^{*}$ achieves a strict local minimum of the height on the set of $n$-dimensional flat tori of covolume $1$. 
 If moreover $\zeta(L,s) <0$ for $0<s<\frac{n}{2}$, then $L$ is $\zeta$-extreme at $s$ for any
$s>0$, $s \neq \frac{n}{2}$.
 
\end{theor}
This theorem applies to $\mathbb{D}_4$, $\mathbb{E}_8$ and
$\Lambda_{24}$, and somehow clarifies Sarnak and Str\"ombergsson's proof. Moreover, it should 
apply to a wider class of lattices, for which the group theoretic
tools used in \cite{Sa-Sa} are not available, but for which one can prove
however that the $4$-design properties hold for all layers. In particular we prove that essentially all the extremal modular lattices (up to certain restrictions on both the dimension and the level) share with $\mathbb{D}_4$, $\mathbb{E}_8$ and
$\Lambda_{24}$ the property of being $\zeta$-extreme at $s$ for any $s> \frac{n}{2}$, $s \neq \frac{n}{2}$ (see
proposition \ref{ml} in section \ref{ex}). 

The paper is organized as follows : in Section \ref{basics} we
collect some preliminary results about lattices, spherical designs and Epstein's
$\zeta$ function. Section \ref{proof} contains
the proof of the main theorem, of which we present several
examples of application in Section \ref{ex}. In Section \ref{mtheta}, we derive a statement about the minima of theta functions, similar to our main theorem for the Epstein zeta function.

\subsection*{Notation}
We denote $x\cdot y$ the usual scalar product of the vectors $x$ and $y$ in $\R^n$, and $\Vert x\Vert$ the associated norm.
We write vectors in $\R^n$ as column vectors. The transpose of a matrix is denoted by the superscript ${}'$. Also, if $A$ is a square $n$-by-$n$ matrix, and $x$ a vector in $\R^n$, the notation $A[x]$ stands for the product $x'Ax$.

The set of $n$-by-$n$ symmetric matrices with real entries is denoted $S_n(\R)$. It is endowed with its canonical scalar product
\begin{equation}
\left\langle A,B\right\rangle :=\Tr AB \ \ \ A,B \in S_n(\R).
\end{equation} 

If $A \in S_n(\R)$ is fixed, one associates with any $x \in \R^n$ such that $A[x] \neq 0$ the $n$-by-$n$ symmetric matrix 
\begin{equation}
\widehat{x} _A:=\dfrac{xx'}{A[x]}.
\end{equation} 
One has $\left\langle H,\widehat{x} _A\right\rangle =\dfrac{H[x]}{A[x]}$ for any $H \in S_n(\R)$.

\section{Basics.}\label{basics}
\subsection{Lattices and quadratic forms.}
Throughout the paper, we denote $\mathcal L_n$ (resp. $\mathcal L_n^{\circ} $) the set of Euclidean lattices (resp. of covolume $1$) in $\R^n$, and $\mathcal P_n$ (resp. $\mathcal P_n^{\circ} $) the cone of positive definite quadratic forms (resp. of determinant $1$) in $n$ variables. We identify a quadratic form in $\mathcal P_n$ with its matrix in the canonical basis of $\R^n$. The map $P \mapsto P'P$ induces a bijection from $O_n(\R)~\setminus~\GL_n(\R)$ onto $\mathcal P_n$ and we thus identify these two sets. Similarly, we identify $\mathcal L_n$ with the quotient $\GL_n(\R)/\GL_n(\Z)$, associating the lattice $L=P\Z^n$ the coset $P \GL_n(\Z) \in \GL_n(\R)/\GL_n(\Z)$. There is a well-known "dictionary" between Euclidean lattices and positive definite quadratic forms which is summarized in the following diagram
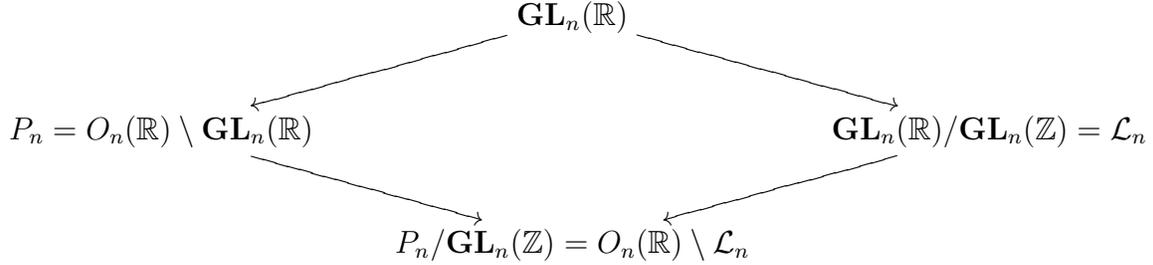
\begin{equation*}
\xymatrix{& \GL_n(\R)\ar[ld] \ar[rd]\\ P_n=O_n(\R) \setminus \GL_n(\R) \ar[rd] && \GL_n(\R)/\GL_n(\Z) =\mathcal L_n \ar[ld]\\&P_n/\GL_n(\Z) =O_n(\R) \setminus \mathcal L_n &}
\end{equation*} 

This enables to formulate every definition and statement in either of these languages but, depending on the context, one point of view is often better than the other. In particular, the proof of the main theorem is easier to write in terms of quadratic forms. 

Let $L=P\Z^n$ be a lattice in $\R^n$, and $A=P'P$ the corresponding quadratic form. We define the sequence $m_1(L) < m_2(L) < \cdots$ of squared lengths of non zero vectors in $L$, ranged in increasing order (this is not the same, in general, as the \textit{successive minima} of $L$). The \textit{$k$-th layer} of $L$, for $k \in \NN \setminus \left\lbrace 0\right\rbrace$, is defined as
\begin{equation}
M_k(L):=\left\lbrace x \in L \mid x\cdot x = m_k(L)\right\rbrace
\end{equation} 
and we set 
\begin{equation}
a_k(L)=\vert M_k(L) \vert.
\end{equation} 
One defines similarly the sequence $\left( m_k(A)\right) _{k \in \NN \setminus \left\lbrace 0\right\rbrace}$ of successive values achieved by the quadratic form $A$ on $\Z^n \setminus \left\lbrace 0\right\rbrace$, as well as the associated layers $M_k(A):=\left\lbrace x \in \Z^n \mid A[x] = m_k(A)\right\rbrace$. 
 
\subsection{Spherical designs.}\label{designs} 
We collect in this section some properties that we will need in the sequel. We take as original definition of a spherical $t$-design the one we gave in the introduction of this text. For sake of completeness, we recall first a proposition to be found in \cite{V}, which provides several characterizations of spherical designs. We recall that a polynomial $P(x_1, \cdots,x_n)$ is harmonic if $\Delta P =0$, where $\Delta$ is the usual laplace operator $\Delta=\sum_{i=1}^n\dfrac{\partial^2}{\partial x_i^2}$.
\begin{proposition} [Venkov {\cite[Th\'eor\`eme 3.2.]{V}}]\label{venkov} 
Let $X$ be a finite subset of $\Sph^{n-1}$ and $t$ an even positive integer. Assume that $X$ is symmetric about $0$, \ie $X=-X$. Then the following properties are equivalent :
\begin{enumerate} 
\item $X$ is a $t$-design.
\item For all non constant harmonic polynomial $P(x)$ of degree $\leq t$, $\sum_{x \in X}f(x)=0$.
\item There exists a constant $c$ such that $\forall \alpha \in \R^n \sum_{x \in X} (x \cdot \alpha) ^t=c (\alpha \cdot \alpha) ^{\frac{t}{2}}$.
\end{enumerate} 
\proof See \cite[Th\'eor\`eme 3.2.]{V}. The only difference is that the condition $\forall \alpha \in \R^n  \sum_{x \in X} (x \cdot \alpha)^i =0$, where $i=t-1$, which should appear in (2), is automatically satisfied here since $X$ is symmetric about the origin. \qed
\end{proposition}
The next proposition is the key to the proof of our main theorem. It is a formulation of the property that all layers of a lattice $L$ hold a $2$- (resp. $4$-) design, in terms of its zeta function. 

\begin{proposition}\label{crit} 
Let $L=P\Z^n$ be a lattice in $\R^n$, and $A=P'P$.
\begin{enumerate} 
\item[1.] The following conditions are equivalent
\begin{enumerate}
\item All layers of $L$ hold a $2$-design.
\item For all $s \in \C$ with $\Re s> \frac{n}{2}$, $\sum_{y \in L \setminus \{0\}}\dfrac{yy'}{\Vert y \Vert ^{2(s+1)}}=\dfrac{\zeta(L,s)}{n}I_n$.
\item For all $s \in \C$ with $\Re s> \frac{n}{2}$, $\sum_{x \in \Z^n \setminus \{0\}} \dfrac{\widehat{x} _A}{A[x]^s}=\dfrac{\zeta(A,s)}{n}A^{-1}$.
\end{enumerate} 

\item[2.] The following conditions are equivalent
\begin{enumerate}
\item All layers of $L$ hold a $4$-design.
\item For all $s \in \C$ with $\Re s> \frac{n}{2}$, for all $H \in S_ n (\R)$, 
\begin{equation*}
\sum_{y \in L \setminus \{0\}} \dfrac{H[y]^2}{\Vert y \Vert ^{2(s+2)}}=\dfrac{\zeta(L,s)}{n(n+2)} ((\Tr H)^2 + 2 \Tr H^2).
\end{equation*} 
\item For all $s \in \C$ with $\Re s> \frac{n}{2}$, for all $H \in S_n (\R) $, 
\begin{equation*}
\sum_{x \in \Z^n \setminus \{0\}} \dfrac{\left\langle H, \widehat{x} _A\right\rangle ^{2}}{A[x]^s}=\dfrac{\zeta(A,s)}{n(n+2)}((\Tr A^{-1}H)^2 + 2 \Tr (A^{-1}H)^2).
\end{equation*}  
\end{enumerate} 
\end{enumerate} 
\end{proposition}
\proof 
1. The equivalence between $(b)$ and $(c)$ is straightforward, using the dictionary lattices/quadratic forms, so it is enough to prove the equivalence between $(a)$ and $(b)$. Assume that all layers of $L$ are $2$-designs. From proposition \ref{venkov}, this means that for all $k \in \NN \setminus \{0\}$, there exists a constant $c_k$ such that
\begin{equation}\label{2d} 
\forall \alpha \in \R^n \sum_{y \in M_k(L)} (y \cdot \alpha) ^2=c_k (\alpha \cdot \alpha) ^2 
\end{equation} 

This equation may be viewed as an equality between quadratic forms in $\alpha$, which, once written matricially, reads
\begin{equation}\label{2dbis} 
\sum_{y \in M_k(L)} yy'=c_k I_n .
 \end{equation} 
It remains to compute the constant $c_k$, which is achieved by taking the trace. This yields 
\begin{equation*}
c_k=\dfrac{a_k(L) m_k(L)}{n}.
\end{equation*} 
Adding up the contributions of all $k \in \NN \setminus \{0\}$, we get the foreseen relation 
\begin{equation}
\sum_{y \in L \setminus \{0\}} \dfrac{yy'}{\Vert y \Vert ^{2(s+1)}}=\dfrac{\zeta(L,s)}{n}I_n.
\end{equation} 
Assume conversely that (b) holds. This is an equality between Dirichlet series, so it must hold coefficientwise. Identifying the coefficients of $m_k(L)^{-s}$ in both sides of (b) yields a relation similar to (\ref{2dbis}), hence a $2$-design relation for the set of vectors of squared length $m_k(L)$.

2. Again, it is enough to prove the equivalence between $(a)$ and $(b)$. Assume first that all layers of $L$ are $4$-designs. For fixed $H \in \mathcal T _A$ we set $P_H(x):=\left\langle H, xx' \right\rangle^2=H[x]^2$, $x \in \R^n$. It is a homogenous polynomial of degree $4$ in $x$. From \cite[Th\'eor\`eme 2.1.]{V}, it decomposes as
\begin{equation}\label{harm} 
P_H(x)=P_4(x)+\Vert x \Vert ^2 P_2(x)+\Vert x \Vert ^4 P_0(x),
\end{equation} 
where $P_i(x)$ is an harmonic polynomial of degree $i$, \ie $\Delta P_i=0$, depending on $H$.
From proposition \ref{venkov}(2), the property that each layer of $L$ holds a $4$-design implies that for $i=2,4$ and for all $k \in \NN \setminus \{0\}$
\begin{equation*}
\sum_{x \in M_k(L)}P_i(x) =0.
\end{equation*}
Consequently, using (\ref{harm}) we obtain
\begin{equation*}
\forall k \in \NN \setminus \{0\}, \sum_{x \in M_k(L)}P_H(x) = \sum_{x \in M_k(L)}\Vert x \Vert ^4 P_0(x)=m_k(L)^2 \# M_k(L) \cdot P_0
\end{equation*} 
since $P_0(x)=P_0$ is a constant. So what we finally need is to compute $P_0$ in (\ref{harm}).
Applying $\Delta$ to the right-hand side of (\ref{harm}) twice, we obtain since the $P_i$ are harmonic, 
\begin{equation}\label{first} 
\Delta^2 P_H(x)= 8n(n+2)P_0.
\end{equation} 
On the other hand, one can compute $\Delta^2 P_H(x)$ directly from the definition of $P_H(x)=H[x]^2$. Easy calculations yield
\begin{equation*}
\Delta H[x] = 2 \Tr H \text{ and } 
\Delta H[x]^2= 4 (\Tr H) H[x] +8 H^2[x]
\end{equation*} 
whence finally
\begin{equation}\label{second} 
\Delta^2 P_H(x)=\Delta ^2 H[x]^2= 8 (\Tr H)^2 + 16 \Tr H^2.
\end{equation} 
Comparing (\ref{first}) and (\ref{second}), we get
\begin{equation}
P_0 = \dfrac{ (\Tr H)^2 + 2 \Tr H^2}{n(n+2)}.
\end{equation} 
Adding up the contributions of all $k \in \NN \setminus \{0\}$ yields formula (b).
Assume conversely that (b) holds. Applying it to $H=\alpha \alpha'$, for a given $\alpha \in \R^n$, we obtain
\begin{equation*}
\sum_{y \in L \setminus \{0\}} \dfrac{(y \cdot \alpha)^4}{\Vert y \Vert ^{2(s+2)}}=\dfrac{\zeta(L,s)}{n(n+2)} 3(\alpha \cdot \alpha)^2
\end{equation*} 
and identifying the coefficient of $m_k(L)^{-s}$ on both sides for all $k \in \NN \setminus \{0\}$, we get
\begin{equation*}
\sum_{y \in M_k(L)} (y \cdot \alpha)^4=\dfrac{3a_k(L) m_k(L)^2}{n(n+2)}(\alpha \cdot \alpha)^2
\end{equation*} 
whence the conclusion.
\qed
 
\begin{remark}
Lattices satisfying any of the equivalent conditions 1(a) and 1(b) in Proposition \ref{crit} are called "strongly critical" in \cite{DR}. This property partly characterizes lattices that are $\zeta$-extreme at $s$ for
any large enough $s$. More precisely, one may reformulate \cite[Theorem 4]{DR} as

\begin{theo}[Delone \& Ryshkov {\cite[Theorem 4]{DR}}]\label{dr}
The following conditions for $L \in \mathcal L_n^{\circ} $ are equivalent.
\begin{enumerate}
\item There exists $s_0>0$ such that $L$ is $\zeta$-extreme at $s$ for
any $s >s_o$.
\item $L$ is perfect, and all layers of $L$ hold a $2$-design.
\end{enumerate} 
\end{theo}
(recall that a lattice $L$ is {\rm perfect} if $\sum_{x\in M_1(L)}\R xx'~=~S_n(\R)$).
\end{remark}
\section{Proof of Theorem \ref{mt}.}\label{proof}
We view $\mathcal P_n^{\circ} $ as a
differentiable submanifold of $S_n(\R)$. The tangent space $\mathcal T _A $ at any point $A$ identifies with the set $\{H \in
S_n(\R) \mid \lb A^{-1}, H \rb = 0\}$. Moreover, the exponential map $H
\mapsto e_A(H)=A\exp (A^{-1}H)$ induces a local diffeomorphism from $\mathcal T _A $ onto $\mathcal P_n^{\circ} $. Consequently we have to study the local behaviour of the map $H \mapsto \zeta(e_A(H),s)$, $H \in \mathcal T _A $ for fixed $s>0$. A simple calculation, based on the Taylor expansion of the exponential, yields 
\begin{eqnarray}\label{taylor} 
\zeta(e_A(H),s)&=& \zeta(A,s)-s\left\langle H, \sum {}' \dfrac{\widehat{x} _A}{A[x]^s}\right\rangle \\
\nonumber & & {}+\frac{s}{2}\left[ (s+1)\sum {}'\dfrac{\left\langle H, \widehat{x} _A\right\rangle ^{2}}{A[x]^s}-\left\langle HA^{-1}H,\sum {}' \dfrac{\widehat{x} _A}{A[x]^s}\right\rangle \right] + o(\Vert H^2 \Vert),
\end{eqnarray}
in which the abbreviated notation $\sum {}'$ stands for $\sum_{x \in \Z^n \setminus \{0\}}$. As explained in the previous section (Proposition \ref{crit}), the property that all layers of $L$ are $2$-designs is equivalent to the relation
\begin{equation}
\sum {}' \dfrac{\widehat{x} _A}{A[x]^s}=\dfrac{\zeta(A,s)}{n}A^{-1}.
\end{equation} 
For $H \in \mathcal T _A$ it implies
\begin{equation}
\left\langle H, \sum {}' \dfrac{\widehat{x} _A}{A[x]^s}\right\rangle =0
\end{equation} 
and 
\begin{equation}
\left\langle HA^{-1}H,\sum {}' \dfrac{\widehat{x} _A}{A[x]^s}\right\rangle =\dfrac{\zeta(A,s)}{n} \Tr(A^{-1}H)^2.
\end{equation}  
Next we use the assumption that all layers of $L$ are $4$-designs to compute the term $\sum {}'\dfrac{\left\langle H, \widehat{x} _A\right\rangle ^{2}}{A[x]^s}$. From Proposition \ref{crit} we have, for $H \in \mathcal T _A$,
\begin{equation}
\sum {}' \dfrac{\left\langle H, \widehat{x} _A\right\rangle ^{2}}{A[x]^s}=\dfrac{\zeta(A,s)}{n(n+2)}((\Tr A^{-1}H)^2 + 2 \Tr (A^{-1}H)^2)= \dfrac{2\zeta(A,s)}{n(n+2)} \Tr (A^{-1}H)^2.
\end{equation}
Inserting the last three formulas into (\ref{taylor}), we obtain
\begin{equation}\label{finaleq} 
\zeta(e_A(H),s) =\zeta(A,s)\left[1+\dfrac{s(s-\frac{n}{2})}{n(n+2)}\Tr (A^{-1}H)^2\right] + o(\Vert H^2 \Vert).
\end{equation} 
Consequently, the assertion that $A$ achieves a strict local minimum on $\mathcal P_n^{\circ} $ of the map $A \mapsto \zeta(A,s)$ is equivalent to the fact that $\zeta(A,s)\dfrac{s(s-\frac{n}{2})}{n(n+2)}>0$. This is clearly the case if $s>\dfrac{n}{2}$, while for $0<s<\frac{n}{2}$ this is equivalent to $\zeta(A,s)<0$. 

As for the assertion on the height function, we just have to differentiate (\ref{finaleq}) with respect to $s$ to get 
\begin{eqnarray} 
\frac{d}{ds}\zeta(e_A(H),s) _{\vert s=0}&=&\zeta'(A,0) +\zeta(A,0) \dfrac{-1}{n+2}\Tr (A^{-1}H)^2+ o(\Vert H^2 \Vert)\\&=&\zeta'(A,0) +\dfrac{1}{n+2}\Tr (A^{-1}H)^2+ o(\Vert H^2 \Vert),
\end{eqnarray} 
since $\zeta(A,0) =-1$, whence the conclusion. This finishes the proof of Theorem \ref{mt}. 
\qed

\section{Examples.}\label{ex}
In order to avoid rescaling systematically all the lattices appearing in the examples below to covolume $1$, we will use the slightly abusive formulation "$L$ is $\zeta$-extreme" to mean that "$L$ \textit{rescaled to covolume }$1$ is $\zeta$-extreme". Similarly, we say that "the torus associated with $L$ achieves a local minimum of the height function" to mean "a local minimum on the set of flat tori \textit{of the same covolume } $\det L$".

Before we give some explicit examples, we wish to give a first comparison between our criterion for a lattice to be $\zeta$-extreme and Sarnak-Str\"ombgergsson's one. For the proof of \cite[Theorem 1]{Sa-Sa}, one considers the space $\Sym^f \Sym^2 (\R^n)$ for $f=0,1,2, \dots$ endowed with the standard action of $O(n)$, and define $f(L)$ to be the largest integer such that $\Sym^f \Sym^2 (\R^n)^{O(n)} = \Sym^f \Sym^2 (\R^n)^{\Aut(L)}$. Then it is proven that if $f(L) \geq 2$ then $L$ is $\zeta$-extreme for $s>\frac{n}{2}$. As noticed by the authors, the determination of $f(L)$ is related to the somewhat more classical problem of determining the largest integer $t(L)$ such that $\Sym^t(\R^n)^{O(n)}=\Sym^t(\R^n)^{\Aut(L)}$, a question which is itself connected with the existence of spherical designs in lattices. To be more precise, one has the following result

\begin{theo}[Goethals and Seidel{\cite[Th\'eor\`eme 6.1.]{GS}}]\label{gs} 
The following conditions for a finite subgroup $G$ of $O(n)$ are equivalent :
\begin{enumerate}
\item $\Sym^t(\R^n)^{O(n)}=\Sym^t(\R^n)^{G}$.
\item Any orbit $G \cdot a$ of a point $a \in \Sph^{n-1}$ is a $t$-design.
\end{enumerate} 
\end{theo}
The combination of this with Theorem \ref{mt} leads to the following corollary
\begin{corollary}\label{cmt} 
If $t(L) \geq 4$ then $L$ is $\zeta$-extreme for $s>\frac{n}{2}$, and the torus associated with $L^{*}$ achieves a local minimum of the height function.
\end{corollary}
\proof For any $k \in \NN \setminus \left\lbrace 0\right\rbrace$ and any $x \in M_k(L)$, the orbit of $x$ under $\Aut(L)$ is a $4$-design, so that $M_k(L)$ is itself a $4$-design, as a union of $4$-designs. \qed

Note that the assumption that $t(L) \geq 4$ is equivalent to $t(L) \geq 5$, since $t(L)$ is easily seen to be odd.  
Sarnak and Str\"ombergsson pointed out that 
\begin{equation}\label{ft} 
f(L) \leq \dfrac{t(L)-1}{2},
\end{equation} 
so that $f(L) \geq 2$ actually implies $t(L) \geq 5$. However, they also observed that (\ref{ft}) is in general a strict inequality, so that the assumption $f(L) \geq 2$ is in general stronger than the assumption of Corollary \ref{cmt}, which is itself stronger than the assumption of Theorem \ref{mt}.

A good illustration of the combination of our criterion with Goethal and Seidel's Theorem, is obtained with the family of Barnes-Wall lattices. Let us briefly recall their definition : if $n=2^k$, we consider an orthonormal basis $\left( e_u \right) _{u \in \F_{2^k}}$ of $\R^n$ indexed by the elements of $F_{2^k}$ and set
\begin{equation}
BW_n:=\left\langle 2^{\lfloor\frac{k-d+1}{2}\rfloor}\sum_{u \in U}e_u\right\rangle _{\Z} \subset \R^n
\end{equation} 
where $U$ runs through the set of affine subspaces of $\R^n$ and $d$ stands for the dimension of $U$. This defines for any $n$ an isodual lattice, \ie isometric to its dual. These lattices are very interesting inasmuch they form one of the very few infinite families of lattices for which explicit computations can be made (density, kissing number, automorphism group etc.), although they do not provide, in dimension $\geq 32$, the best known lattice packings. Various explicit descriptions of the automorphism group of $BW_n$ are known (see for instance \cite{BE}), and its polynomial invariants are computed in \cite{Ba}. Altogether, this leads to the following proposition :
\begin{proposition}
For $k\geq 3$, the Barnes-Wall lattice $BW_{2^k}$ is $\zeta$-extreme for any $s>\frac{n}{2}$ and the associated torus achieves a strict local minimum of the height function.
\end{proposition}
\proof From  \cite[Corollary 5.1.]{Ba}, we see that $t(BW_{2^k}) \geq 6$, whence the conclusion using corollary \ref{cmt} (notice that for an isodual lattice $L$, the tori associated with $L$ or its dual are the same, up to scaling). \qed

Other examples of lattices $L$ satisfying $t(L) \geq 4$ (and thus our main theorem) are the root lattices $\mathbb E_6$ and $\mathbb E_7$, as well are their duals. Finally, a list of lattices in dimension $\leq 26$ to which this argument applies is provided in \cite[table 1]{Ba} (note that only the lattices pertaining to what is called case (1) and (3) there are suitable).

In the examples below, we prove that some lattices $L$ are $\zeta$-extreme without computing neither $f(L)$ nor $t(L)$. Instead, we refer to the paper \cite{BV} by Bachoc and Venkov, where the existence of spherical designs in certain lattices is proven using modular forms. The lattices dealt with in that paper 
pertain to Quebbemann's theory of modular lattices (see \cite{Q} or \cite{Sch-Sch}), from which we recall the main definitions and results. A lattice $L$ in $\mathcal L_n$ is $\ell$-\textit{modular} ($\ell>0$) if
\begin{enumerate}
\item [(i)] $L$ is \textit{even}, \ie $x\cdot x \in 2\Z$ for all $x \in L$.
\item [(ii)] $L$ is isometric to $\sqrt{\ell}L^{*}$.
\end{enumerate} 
Assume furthermore that $\ell \in \left\lbrace 1,2,3,5,7,11,23\right\rbrace$. Then ( \cite[Theorem 2.1.]{Sch-Sch}) the minimum of an $\ell$-modular lattice $L$ satisfies
\begin{equation}\label{extr} 
\min L \leq 2(1+\lfloor \dfrac{n(1+\ell)}{48}\rfloor)
\end{equation} 
An $\ell$-modular lattice $L$ for which equality holds in (\ref{extr}) is called \textit{extremal}. The theta series, and more generally the theta series with spherical coefficients of $\ell$-modular lattices belong to a a certain algebra of modular forms, which can be described explicitely when $\ell$ belongs to the set above. From this description, Bachoc and Venkov deduce various results about the existence of spherical designs in extremal modular lattices (see \cite[Corollary 3.1.]{BV}). Applying their result together with Theorem \ref{mt}, one easily derives the following proposition :
\begin{proposition}\label{ml} 
Let $L$ be an extremal $\ell$-modular lattice of dimension $n$ 
such that $\ell=1$ and $n \equiv 0,8 \mod 24$, or $\ell=2$ and $n \equiv 0,4 \mod 16$, or $\ell=3$ and $n \equiv 0,2 \mod 12$. Then $L$ is $\zeta$-extreme for any $s>\frac{n}{2}$ and the associated torus achieves a strict local minimum of the height function.
\end{proposition}
\proof From \cite[Corollary 3.1.]{BV}, all layers of such a lattice hold a $4$-design, whence the conclusion.  The proof of this fact uses theta series with spherical coefficients. \qed

The previous proposition applies to the $\mathbb{E}_8$ and Leech lattices ($l=1$), as well as the $\mathbb{D}_4$ lattice ($l=2$), recovering Sarnak and Str\"ombergsson's result, at least for $s>\frac{n}{2}$. But this also applies for instance to the hexagonal lattice $\mathbb{A}_2$, to all extremal even unimodular lattices in dimension $32$ and $48$, to the Coxeter-Todd lattice $\mathrm{K}_{12}$ or the Barnes-Wall lattice $\Lambda_{16}$ to cite a few. Note that the occurrence of the hexagonal lattice $\mathbb{A}_2$ in this list is not a surprise, since it is known to achieve a \textit{global} minimum of the map $L \mapsto \zeta(L,s)$ at $s$ for any $s>0$, $s \neq \frac{n}{2}$ (see \cite{Di}, \cite{En} and \cite{Ra}).
\section{Minima of theta functions.}\label{mtheta} 
In this section, we investigate the question of the minima of theta functions, which is closely related to the subject dealt with so far. Recall that the theta function of a lattice $L$ is defined as  
\begin{equation}\label{theta} 
\Theta_L(z)=\sum _{l \in L}e^{\pi i z \Vert l \vert^2}, \ \text{ for } z \in \C, \ \Im z >0.
\end{equation} 
The theta an zeta functions of a lattice are related through Mellin transform, namely one has, for $s \in \C$ with $\Re s >\frac{n}{2}$,
\begin{equation}\label{mellin} 
\Gamma(s) \pi^{-s}\zeta(L,s)=\mathcal M (\Theta_L(iy)-1):=\int_{0}^{+\infty}(\Theta_L(iy)-1)y^{s-1}dy.
\end{equation} 

For fixed $y>0$, we ask for lattices in $\mathcal L_n^{\circ}$ minimizing $\Theta_L(iy)$. Sarnak and Str\"ombergsson proved \cite[Proposition 2]{Sa-Sa} that for any $y>0$, the $\mathbb D_4$ lattice (rescaled to covolume $1$),  the $\mathbb{E}_8$ lattice and the Leech lattice achieve a strict local minimum of $\Theta_L(iy)$. Following the same line as in the proof of our main theorem, we can prove the following result :
\begin{proposition}\label{minitheta}  Let $L_0 \in \mathcal L_n ^{\circ}$ be such that all its \textrm{layers}
hold a $4$-design. Then, for any fixed $y>\frac{\frac{n}{2}+1}{\pi m_1(L_0)}$, the map $L \mapsto \Theta(L,iy)$, $L \in \mathcal L_n ^{\circ}$, has a strict local minimum at $L_0$.
\end{proposition}
\proof As before, we give the proof in terms of positive definite quadratic forms. If $B$ is the positive definite symmetric matrix, defined up to $\GL_n(\Z)$ equivalence, corresponding to a lattice $L$, one defines $\Theta_B(z)=\Theta_L(z)=\sum_{x \in \Z^n}e^{\pi i z B[x]}$, for $z \in \C$ with $\Im z >0$. Letting $A$ be the positive definite symmetric matrix associated to $L_0$, we parametrize locally the set $\mathcal P_n^{\circ}$ via the exponential map $e_A$ as in the proof of the main theorem, and we are led to study the local behaviour of the map $H \mapsto \Theta_{e_A(H)}(iy)$, $H \in \mathcal T _A$. Under the conditions of the proposition, equation (\ref{finaleq}) holds. Applying inverse Mellin transform to this equation, we get 
\begin{equation}\label{imellin} 
\Theta_{e_A(H)}(iy)=\Theta_{A}(iy)+\dfrac{\Tr (A^{-1}H)^2 }{n(n+2)}\mathcal M^{-1}\left( s(s-\frac{n}{2}) \Gamma(s) \pi^{-s}\zeta(A,s)\right)  + o(\Vert H^2 \Vert).
\end{equation} 
Using elementary properties of the Mellin inverse transform, we thus find
\begin{equation}
\Theta_{e_A(H)}(iy)=\Theta_{A}(iy)+\dfrac{\Tr (A^{-1}H)^2  }{n(n+2)}\sum_{x \in \Z^n}y\pi A[x]\left(y\pi A[x]-(\frac{n}{2}+1) \right)  e^{-\pi y A[x]}+ o(\Vert H^2 \Vert).
\end{equation} 
In order to conclude, it is enough to show that the sum $\sum_{x \in \Z^n}y\pi A[x]\left(y\pi A[x]-(\frac{n}{2}+1) \right)  e^{-\pi y A[x]}$ is positive, which is obviuosly the case if $y>\frac{\frac{n}{2}+1}{\pi m_1(L_0)}$, since then each term of the sum is positive. \qed
\begin{remark}
This proposition applies to all the examples dealt with in the previous section. A more careful analysis of the sign of the sum $\sum_{x \in \Z^n}y\pi A[x]\left(y\pi A[x]-(\frac{n}{2}+1) \right)  e^{-\pi y A[x]}$ would allow to extend the range (in $y$) of validity of the proposition, as done by Sarnak and Str\"ombergsson in the case of $\mathbb D_4$, $\mathbb{E}_8$ and Leech lattice. 
 
\end{remark}
\section{Final remarks.}\label{fr}
We conclude with some remarks and open questions. 
\begin{enumerate}
\item In the examples quoted above, we applied our main theorem to derive $\zeta$-extremality at $s$ for any $s>\frac{n}{2}$. To get the same result for $0<s<\frac{n}{2}$, one has to prove that $\zeta(L,s) <0$  in that range, which was done by  Sarnak and Str\"ombergsson in the case of $\mathbb{D}_4$, $\mathbb{E}_8$ and $\Lambda_{24}$. Unfortunately, we don't know of any 'uniform' way to prove this property, so that a case-by-case proof would be necessary to deal with the examples of the previous section. However, as pointed out to me by P. Sarnak, it can be shown, using an argument due to A. Terras \cite[Theorem 1]{T} that the Epstein zeta function of the Barnes-Wall lattices $BW_n$ do have a zero in $(0,\frac{n}{2})$ for large enough $n$, so that the extremality does not hold for all $0<s<\frac{n}{2}$. The same argument would also apply to extremal modular lattices of large enough dimension, provided that they exist (for fixed level $\ell$, the dimension of an hypothetical extremal modular lattice is bounded, see \cite[Theorem 2.1. (ii)]{Sch-Sch}).
\item The condition that all the layers of a given lattice
hold a $4$-design is rather strong. Lattices for which the first layer (minimal vectors) hold a $4$-design, the so-called 'strongly perfect lattices', have been classified in dimensions up to $12$ (see \cite{V}, \cite{NV1}, \cite{NV2}). In turns out that in dimension $3$, $5$, and $9$, for instance, such a lattice  (and \textit{a fortiori} a lattice all the layers of which hold a $4$-design) does not exist. As for the weaker condition that all the layers of a given lattice
hold a $2$-design, which is necessary for the lattice to be $\zeta$-extremal at $s$ for all big enough $s$, according to Delone and Ryshkov's theorem, it is easy to find examples in any dimension : for instance, all irreducible root lattices hold this property. However, it is still not clear that a lattice achieving a global minimum of the function $L \mapsto \zeta(L,s)$ for any $s > \dfrac{n}{2}$ (or even for all large enough $s$) should exist.
\item The situation for the height function is perhaps more intriguing. Indeed, it is known (see \cite{Ch}) that in a given dimension $n$, the height function, restricted to flat tori, achieves a global minimum. On the other hand, in dimension $3$, $5$, and $9$ for instance, there is no hope to find this minimum using the criterion of Theorem \ref{mt}. Consequently, the \textit{right} characterization of local minima of the height function is still to be found (our condition is too strong). Recall that such a characterization for the local maxima of the density of lattice-sphere packings is known, due to Vorono\"i (see \cite{Vor} or \cite[Chapter 3]{M}).
\end{enumerate} 
\section*{Acknowledgements.} I would like to thank Peter Sarnak for his comments on a preliminary version of this work, which led in particular to the statement of Proposition \ref{minitheta}.


\newpage
\begin{thebibliography}{11}
 
\bibitem{Ba} {C. Bachoc},  
Designs, groups and lattices.  J. Th\'eor. Nombres Bordeaux {\bf 17} (2005), no. 1, 25--44.  

\bibitem{BV} {C. Bachoc, B. Venkov}, Modular forms, lattices and spherical designs. {\em R\'eseaux euclidiens, designs sph\'eriques et formes modulaires}, 10--86, Monogr. Enseign. Math., 37, Enseignement Math., Geneva, 2001.

\bibitem{BE}{M. Brou\'e, M. Enguehard}, Une famille infinie de formes quadratiques enti\`eres; leurs groupes d'automorphismes., Ann. Sci. \'Ecole Norm. Sup. (4)  \textbf{6} (1973), 17--51.

\bibitem{Ca} {J.W.S. Cassels}, On a problem of Rankin about the
 Epstein zeta function, Proc. Glasgow Math. Assoc. {\bf 4} (1959),
 73--80, {\bf 6} (1963), 116.
 
 \bibitem{Ch}{P. Chiu}, Height of flat tori.  Proc. Amer. Math. Soc.  {\bf 125}  (1997),  no. 3, 723--730.

\bibitem{DR} { B. N. Delone, S. S. Ry\v skov}, A contribution to the theory of the extrema of a multi-dimensional $\zeta$-function. Dokl. Akad. Nauk SSSR 173 991--994 (Russian); translated as Soviet Math. Dokl. 8 1967 499--503.

\bibitem{Di} {P.H. Diananda}, Notes on two lemmas concerning the
 Epstein-zeta function, Proc. Glasgow Math. Assoc. {\bf 6} (1964),
 202--204.

\bibitem{En} {V. Ennola}, A lemma about the
 Epstein-zeta function, Proc. Glasgow Math. Assoc. {\bf 6} (1964),
 198--201.
 
\bibitem{GS} {J.-M. Goethals, J.J. Seidel}, Spherical designs. Proc. Sympos. Pure Math. {\bf 34} (1979), 255--272. 

\bibitem{M}{J. Martinet}, {\em Perfect lattices in Euclidean spaces.} Grundlehren der Mathematischen Wissenschaften, 327. Springer-Verlag, Berlin, 2003.
 
\bibitem{NV1}{G. Nebe, B. Venkov}, The strongly perfect lattices of dimension 10. Colloque International de Théorie des Nombres (Talence, 1999).  J. Théor. Nombres Bordeaux  {\bf 12}  (2000),  no. 2, 503--518.

\bibitem{NV2}{G. Nebe, B. Venkov}, Low-dimensional strongly perfect lattices. I. The 12-dimensional case.  Enseign. Math. (2)  51  (2005),  no. 1-2, 129--163.
 
\bibitem{Q}{H.G. Quebbemann}, Modular lattices in Euclidean spaces. J. Number Theory 54 (1995), no. 2, 190--202.

\bibitem{Ra} {R.A. Rankin}, A minimum problem for the Epstein zeta function, Proc. Glasgow Math. Assoc. {\bf 1} (1953),
 149--158.
 
\bibitem{Sa} {P. Sarnak}, Determinants of Laplacians; heights and finiteness. {\em Analysis, et cetera},  601--622, Academic Press, Boston, MA, 1990.

\bibitem{Sa-Sa} {P. Sarnak, A. Str\"ombergsson}, Minima of Epstein's Zeta Function and Heights of Flat Tori , Invent. Math. 165 (2006), 115--151.

\bibitem{Sch-Sch}{R. Scharlau, R. Schulze-Pillot}, Extremal lattices. Algorithmic algebra and number theory (Heidelberg, 1997), 139--170, Springer, Berlin, 1999.

\bibitem{So} { S. L. Sobolev }, Formulas for mechanical cubatures in $n$-dimensional space, Dokl. Akad. Nauk SSSR 137 (1961) 527--530.

\bibitem{T} {A. Terras}, The minima of quadratic forms and the behavior of Epstein and Dedekind zeta functions,  J. Number Theory  12  (1980), no. 2, 258--272.

\bibitem{V} {B. Venkov}, R\'eseaux et designs sph\'eriques. {\em R\'eseaux euclidiens, designs sph\'eriques et formes modulaires}, 10--86, Monogr. Enseign. Math., 37, Enseignement Math., Geneva, 2001.

\bibitem{Vor} {G. Vorono\"i}, Nouvelles applications des param\`etres continus \`a la th\'eorie des formes quadratiques : 1 Sur quelques propri\'et\'es des des formes quadratiques parfaites, J. Reine angew. Math. {\bf 133} (1908), 97-178.
 
\end{thebibliography}
\end{document}